\documentclass[12pt,a4paper,final]{amsart}
\allowdisplaybreaks[4]

\usepackage{amsmath}
\usepackage{amssymb}
\usepackage{amsthm}
\usepackage{color}
\usepackage{enumerate}
\usepackage[dvips]{graphicx}
\usepackage{comment}
\usepackage{cleveref}

\newcommand{\T}{\mathbb{T}}

\DeclareMathOperator\SC{SC}

\def\R{\mathbb{R}}

\def\N{\mathbb{N}}

\def\cF{\mathcal{F}}

\def\cA{\mathcal{A}}
\def\bye{\end{document}}
\def\by{\end{proof}\bye}

\def\hello{\begin{document}}
\def\fr{\frac}
\def\disp{\displaystyle}
\def\ga{\alpha}
\def\go{\omega}
\def\gep{\varepsilon}
\def\ep{\gep}
\def\mid{\,:\,}
\def\gb{\beta}
\def\gam{\gamma}
\def\gd{\delta}
\def\gz{\zeta}
\def\gth{\theta}
\def\gk{\kappa}
\def\gl{\lambda}
\def\gL{\Lambda}
\def\gs{\sigma}
\def\gf{\varphi}
\def\tim{\times}
\def\aln{&\,}
\def\ol{\overline}
\def\ul{\underline}
\def\pl{\partial}
\def\hb{\text}
\def\cF{\mathcal{F}}
\def\Int{\mathop{\text{int}}}

\def\gG{\varGamma}
\def\lan{\langle}
\def\ran{\rangle}
\def\cD{\mathcal{D}}
\def\cB{\mathcal{B}}
\def\bcases{\begin{cases}}
\def\ecases{\end{cases}}
\def\balns{\begin{align*}}
\def\ealns{\end{align*}}
\def\balnd{\begin{aligned}}
\def\ealnd{\end{aligned}}
\def\bgat{\begin{gathered}}
\def\egat{\end{gathered}}
\def\1{\mathbf{1}}
\def\bproof{\begin{proof}}

\def\eproof{\end{proof}}

\theoremstyle{definition}
\newtheorem{definition}{Definition}%[section]
\theoremstyle{plain}
\newtheorem{theorem}[definition]{Theorem}
\newtheorem{corollary}[definition]{Corollary}
\newtheorem{lemma}[definition]{Lemma}
\newtheorem{proposition}[definition]{Proposition}
\theoremstyle{remark}
\newtheorem{remark}[definition]{Remark}
\newtheorem{notation}[definition]{Notation}

\def\hr#1{\rule{0pt}{#1 pt}}

\def\red#1{\textcolor{red}{#1}}
\def\blu#1{\textcolor{blue}{#1}}
\def\beq{\begin{equation}}
\def\eeq{\end{equation}}
\def\bthm{\begin{theorem}}
\def\ethm{\end{theorem}}
\def\bproof{\begin{proof}}
\def\eproof{\end{proof}}

\usepackage[abbrev]{amsrefs}
\newcommand\coolrightbrace[2]{%
\left.\vphantom{\begin{matrix} #1 \end{matrix}}\right\}#2}
\def\eqr#1{\eqref{#1}}
\def\bmat{\begin{pmatrix}}
\def\emat{\end{pmatrix}}

\newcommand{\Pmo}{\mathcal{P}^-_{1}}
\newcommand{\Ppo}{\mathcal{P}^+_{1}}
\newcommand{\Pmk}{\mathcal{P}^-_{k}}
\newcommand{\Ppk}{\mathcal{P}^+_{k}}
\def\diag{\operatorname{diag}}
\def\rT{\mathrm{T}}
\newcommand{\Rn}{{\mathbb R}^N}
\def\IN{\text{ in } }

\def\AND{\text{ and }}
\def\FOR{\text{ for }}
\def\FORALL{\text{ for all }}
\def\ON{\text{ on }}
\def\IF{\hb{ if }}
\def\WITH{\text{ with }}

\def\rmb{\mathrm{b}}
\def\I{\mathbb{I}}
\def\du#1{\left\lan#1\right\ran}
\def\bald{\begin{aligned}}
\def\eald{\end{aligned}}
\def\stm{\setminus}
\def\t{\tau}
\def\SC-{\operatorname{SC}^-}
\def\sc-{\operatorname{sc}^-}
\def\lip{\operatorname{lip}}
\def\B{\operatorname{\mathbb{B}}}

\def\Sp{\operatorname{Sp}}
\def\P{\operatorname{\mathbb{P}}}
\def\erf{\eqref}
\def\cG{\mathcal{G}}
\def\cE{\mathcal{E}}
\def\pr{\,^\prime}
\def\gX{\Xi}
\def\gS{\Sigma}
\def\cV{\mathcal{V}}
\def\cW{\mathcal{W}}
\def\Hall{(H_i)_{i\in\I}}\def\Lall{(L_i)_{i\in\I}}
\def\prop{proposition}
\def\thm{theorem}
\def\lem{lemma}

\def\fC{\frak{C}}
\def\fM{\frak{M}}
\def\0{\mathbf{0}}

\renewcommand{\subjclassname}{%
\textup{2010} Mathematics Subject Classification}

\title[Vanishing discount problem]{An example in the vanishing discount problem for 
monotone systems of Hamilton-Jacobi equations}

\author[H. Ishii]{Hitoshi Ishii}
\address[\textsc{Hitoshi Ishii}]{Institute for Mathematics and Computer Science\newline
\indent Tsuda University  \newline
 \indent   2-1-1 Tsuda, Kodaira, Tokyo, 187-8577 Japan.
}
\email{hitoshi.ishii@waseda.jp}

\keywords{systems of Hamilton-Jacobi equations, vanishing discount, full convergence}
\subjclass[2010]{
35B40, %Asymptotic behavior of solutions,
35D40, %Viscosity solutions 
35F50, %Nonlinear first-order systems
49L25 %Viscosity solutions
}

\def\alert#1{\begin{color}{red}#1 \end{color}}

\begin{document}
\maketitle
\begin{abstract} 
In recent years, there have been many contributions to the vanishing discount 
problem for Hamilton-Jacobi equations. 
In the case of the scalar equation, B. Ziliotto [Convergence of the solutions of the discounted Hamilton-Jacobi equation: a counterexample. J. Math. Pures Appl. (9) 128 (2019), 330-338] has given an example of the Hamilton-Jacobi equation having non-convex Hamiltonian in the gradient variable, for which the full convergence of the solutions does not hold as the discount factor tends to zero. 
We give an example of the nonlinear monotone system of Hamilton-Jacobi equations having convex  Hamiltonians in the gradient variable, for which the whole family convergence of the solutions does not hold. 
\end{abstract}

\tableofcontents

\def\rmc{\mathrm{c}}

\section{Introduction}
We consider the system of Hamilton-Jacobi equations
\beq \label{eq:1}
\bcases
\gl u_1(x)+H_1(Du_1(x))+B_1(u_1(x),u_2(x))=0 \ & \IN \T^n, \\
\gl u_2(x)+H_2(Du_2(x))+B_2(u_1(x),u_2(x))=0 \ & \IN \T^n,
\ecases
\eeq
where $\gl>0$ is a given constant, called the discount 
factor, and the functions $H_i : \R^n\to \R$ and $B_i : \R^2 \to R$, with 
$i=1,2$, are given continuous functions.

In a recent paper \cite{IJ}, the authors have investigated the vanishing discount problem for a nonlinear monotone system of Hamilton-Jacobi equations 
\beq\label{eq:1.2}\bcases
\gl u_1(x)+G_1(x,Du_1(x),u_1(x),u_2(x),\ldots,u_m(x))=0 \ & \IN \T^n, \\
\phantom{\gl u_1(x)+G_1(x,Du_1(x),u_1(x),u_2(x)}\vdots &\\
\gl u_m(x)+G_m(x,Du_m(x), u_1(x),u_2(x),\ldots,u_m(x))=0 \ & \IN \T^n,
\ecases
\eeq
and established under some hypotheses on $G_i 
\in C(\T^n\tim\R^n\tim \R^m)$ that, when $u^\gl=(u_1^\gl,\ldots,u_m^\gl)\in C(\T^n)^m$ denoting the (viscosity) solution 
of \erf{eq:1.2}, the whole family $(u^\gl)_{\gl>0}$ converges in $C(\T^n)^m$ 
to some $u^0\in C(\T^n)^m$ as $\gl\to 0+$.  The constant $\gl>0$ in the above system is the so-called discount factor. 

The hypotheses on $G_i$ are the convexity, coercivity, and monotonicity of $G_i$ 
as well as the solvability of \erf{eq:1.2}, with $\gl=0$.  Here the convexity 
of $G_i$ is meant that the functions $\R^n\tim\R^m\ni (p,u)\mapsto 
G_i(x,p,u)$ are convex. We refer to \cite{IJ} for the precise 
statements of the hypotheses.  

Before \cite{IJ}, there have been many contributions to the question about 
the full convergence in the vanishing discount problem, which we refer to 
\cites{IJ, DZ2, DFIZ, IMT1, IMT2, IS, CCIZ} and the references therein. 

In the case of the scalar equation, B. Ziliotto \cite{Zi} has recently shown an example of the Hamilton-Jacobi equation having non-convex Hamiltonian in the gradient variable for which the whole family convergence does not hold.

Our purpose in this paper is that, by adapting the idea of Ziliotto \cite{Zi} 
to the system \erf{eq:1.2}, we give an example of $B_i$ such that, if  $H_1(0)=H_2(0)=0$, then the solutions of 
the system \erf{eq:1.2}, with $\gl>0$, are bounded but do not converge 
to a single point as $\gl\to 0+$.

Motivated by \cite{Zi}, we fix $d>1$, set 
\[
\gamma = \fr{d}{1+d}, 
\]
choose $K\in\N$ so that
\[
\gamma+4^{-K}<1, 
\]
and write  
\begin{gather*}
\I=\{1,\,2\},\qquad 
\N_K=\{k\in\N\mid k\geq K\}, \qquad 
\cA=\{\gamma +4^{-k}\mid k\in \N_K\}\cup\{\gamma\},\qquad
\cB=\{0,\,1\}, 
\\ \intertext{and set for $\ga,\gb\in\R$, }
\bald
C(\ga,\gb)&\,=\left(c_{ij}(\ga,\gb)\right)_{i,j\in\I}
\\&\,:= \bmat \ga+\gb-2\ga\gb & -(\ga+\gb-2\ga\gb)\\
-( \ga+\gb-2\ga\gb)&  \ga+\gb-2\ga\gb  \emat.  
\eald
\end{gather*}
Moreover, we define $L_i : \R^2\to\R$, with $i\in\I$,  by 
\[
L_1(\ga,\gb)=\ga\gb+d^2(1-\ga)(1-\gb),\quad 
L_2(\ga,\gb)=-\ga\gb-d^2(1-\ga)(1-\gb).
\]

Note that $\cA\subset (0,\,1)$ and $\cA, \,\cB$ are compact subsets of $\R$
and that for all $(\ga,\gb)\in \R\tim\R$,
\beq\label{eq:2}\left\{\bald
&c_{11}(\ga,\gb)+c_{12}(\ga,\gb)=c_{21}(\ga,\gb)+c_{22}(\ga,\gb)=0, 
\\ & c_{11}(\ga,\gb)=c_{22}(\ga,\gb)>0 \ \ \IF (\ga,\gb)\in(0,\,1)\tim[0,\,1]. 
\eald \right. \eeq

Let $\gl\in(0,\,1)$. Consider the problem of finding $u=(u_1,u_2)\in\R^2$ such that 
\beq \label{eq:3}
\bcases 
\gl u_1+B_1(u_1,u_2)=0, 
& \\
\gl u_2+B_2(u_1,u_2)=0.&
\ecases
\eeq
where $B_i:\R^2\to \R$, with $i\in\I$, are the continuous functions given by 
\beq\label{eq:3.1}
B_i(u_1,u_2)= 
\max_{\ga\in\cA}\min_{\gb\in\cB}(c_{i1}(\ga,\gb)u_1+c_{i2}(\ga,\gb)u_2-L_i(\ga,\gb)).
\eeq
For later convenience, we set for $(i,\ga,\gb)\in\I\tim [0,\,1]\tim [0,\,1]$ 
and $(u_1,u_2)\in\R^2$, 
\[
b_i(\ga,\gb,u_1,u_2)=c_{i1}(\ga,\gb)u_1+c_{i2}(\ga,\gb)u_2-L_i(\ga,\gb). 
\]
It is clear that the functions $b_i$ are continuous on $[0,\,1]\tim[0,\,1]\tim\R\tim\R$. 

Although our main concern is the system \erf{eq:3}, for the argument below 
we need to treat a more general form of \erf{eq:3}, that is, the system
\beq \label{eq:3.2}\bcases
\gl u_1+A_1(u_1,u_2)=g_1, &\\
\gl u_2+A_2(u_1,u_2)=g_2, &
\ecases
\eeq
where $(g_1,g_2)\in\R^2$ is a given vector, $A_i$, with $i\in\I$, are defined by 
\[
A_i(u_1,u_2)=\max_{\ga\in\cA_i}\min_{\gb\in \cB_i} b_i(\ga,\gb,u_1,u_2),
\]
and $\cA_i, \cB_i$, with $i\in\I$, are given compact subsets of $[0,\,1]$, 
If we take $g_i=0$, $\cA_i=\cA$, and $\cB_i=\cB$ for all $i\in\I$, 
\erf{eq:3.2} is exactly the system \erf{eq:3}.

We note that the mapping $A : (u_1,u_2)\mapsto (A_1(u_1,u_2),A_2(u_1,u_2))$ is monotone in the sense 
that for any $(u_1,u_2),(v_1,v_2)\in\R^2$, if $u_j-v_j\geq u_k-v_k$, where $j\not=k$, 
then $A_j(u_1,u_2)\geq A_j(v_1,v_2)$. Indeed, assuming, for instance, that $u_1-v_1\geq u_2-v_2$,  we observe that for any $\ga\in\cA_1$ and some $\bar\gb(\ga)\in\cB_1$, 
\[
A_1(u_1,u_2)\geq \min_{\gb\in\cB}b_1(\ga,\gb,u_1,u_2) 
= b_1(\ga,\bar\gb(\ga),u_1,u_2),
\]
while for some $\bar \ga\in\cA_1$, 
\[
A_1(v_1,v_2)=\min_{\gb\in\cB}b_1(\bar \ga,\gb,v_1,v_2)
\leq b_1(\bar\ga,\bar\gb(\bar\ga),v_1,v_2).
\]
Combining the above two, using \erf{eq:2}, and writing $\bar\gb=\bar\gb(\bar\ga)$, we deduce that 
\[\bald
A_1(u_1,u_2)-A_1(v_1,v_2)&\,
\geq b_1(\bar\ga,\bar\gb,u_1,u_2)-b_1(\bar\ga,\bar\gb,v_1,v_2)
\\&\,=c_{11}(\bar\ga,\bar\gb)(u_1-v_1)+c_{12}(\bar\ga,\bar\gb)(u_2-v_2)
\\&\,\geq c_{11}(\bar\ga,\bar\gb)(u_1-v_1)+c_{12}(\bar\ga,\bar\gb)(u_1-v_1)=0. 
\eald
\]
This shows that $A$ is monotone.

With the fact that $A$ is monotone, the following proposition is well-known.

\begin{proposition} \label{prop1} There exists a unique solution $(u_1,u_2)\in\R^2$ 
of \erf{eq:3.2}. Moreover, if $(v_1,v_2)\in\R^2$ satisfies 
\begin{align} \label{eq:4}
&\gl v_1+A_1(v_1,v_2)\leq g_1 \ \ \text{ and } \ \  \gl v_2+A_2(v_1,v_2)\leq g_2,
\intertext{(resp., }  
&\gl v_1+A_1(v_1,v_2)\geq g_1
\ \ \text{ and } \ \  \gl v_2+A_2(v_1,v_2)\geq g_2\ ),  \label{eq:5}
\end{align}
then $\,v_1\leq u_1\,$ and $\,v_2\leq u_2\,$ (resp., $\,u_1\leq v_1\,$ and $\,u_2\leq v_2$). 
\end{proposition}

For completeness, we provide below proof of the above proposition.

\bproof We first show the second claim in Proposition \ref{prop1}, that is, 
the comparison claim. Let $(u_1,u_2)\in\R^2$ be a solution of \erf{eq:3.2} 
and let $(v_1,v_2)\in\R^2$ satisfy either \erf{eq:4} or \erf{eq:5}. 

We treat only the case of \erf{eq:4}; the proof of the 
other case is similar. 
 We argue by contradiction, and hence suppose that either 
$v_1>u_1$ or $v_2>u_2$. Assume, for instance, that 
$v_1-u_1\geq v_2-u_2$. By the monotonicity of $A$, we have 
$A_1(v_1,v_2)\geq A_1(u_1,u_2)$. Hence, by \erf{eq:3.2} and \erf{eq:4}, we deduce that 
\[
0\geq \gl (v_1-u_1)+A_1(v_1,v_2)-A_1(u_1,u_2)\geq \gl(v_1-u_1),
\]
which contradicts that $v_1>u_1$.  Thus, we have $v_1\leq u_1$ and $v_2\leq u_2$. 
The uniqueness of the solution of \erf{eq:3.2} is now evident.

Next, we show that there exists a solution of \erf{eq:3.2}. 
We choose a constant $R>0$ such that $\max_{i\in\I}(\|L_i\|_\infty+|g_i|)\leq \gl R$, 
where $\|\cdot\|_\infty$ denotes the sup-norm, and deduce with the help of \erf{eq:2} that $(v_1,v_2)=-(R,R)$ (resp., $(v_1,v_2)=(R,R)$) satisfies \erf{eq:4} 
(resp., \erf{eq:5}).  It is easily seen that $A : \R^2\to\R^2$ is Lipschitz continuous. 
Let $M\geq 0$ be a Lipschitz constant of $A$, and consider the problem of finding 
$(u_1,u_2)\in\R^2$ that satisfies
\beq\label{eq:6}
(\gl+M) u_1+A_1(u_1,u_2)=f_1 \ \ \AND \ \  (\gl+M)
 u_2+A_2(u_1,u_2)=f_2,
\eeq
where $(f_1,f_2)\in\R^2$ is a fixed vector. By the Banach fixed-point theorem, 
\erf{eq:6} has a unique solution. 

We define inductively a sequence of points $(u_1^{(j)},u_2^{(j)})\in\R^2$, with $j\in\N$, by setting 
$(u_1^{(1)},u_2^{(1)})=-(R,R)$ and, when $(u_1^{(j-1)},u_2^{(j-1)})$ is given, 
solving \erf{eq:6} for $(u_1,u_2)$, with 
\[
(f_1,f_2)=(g_1,g_2)+M(u_1^{(j-1)},u_2^{(j-1)}), 
\]
to set $(u_1^{(j)},u_2^{(j)}):=(u_1,u_2)$. 

Note that $(v_1,v_2)=(R,R)$ (resp., $(v_1,v_2)=-(R,R)$) satisfies 
\[\bgat
(\gl+M) v_1+A_1(v_1,v_2)\geq g_1+MR \ \ \AND \ \ (\gl+M)
 v_2+A_2(v_1,v_2)\geq g_2+MR,
\\ (\text{ resp., }\quad (\gl+M) v_1+A_1(v_1,v_2)\leq g_1-MR \ \ \AND \ \ (\gl+M)
 v_2+A_2(v_1,v_2)\leq g_2-MR.\ )
\egat\]
Applying the comparison assertion of Proposition \ref{prop1},  
with $\gl$ and $g_i$ replaced respectively 
by $\gl+M$ and $g_i+Mu_i^{(1)}$, we obtain $-R\leq u_1^{(2)}\leq R$ and $-R\leq u_2^{(2)}\leq R$. 
It is easily seen by induction that $-R\leq u_1^{(j)}\leq R$ 
and $-R\leq u_2^{(j)}\leq R\,$
for all $j\in\N$.  

Moreover, if we assume that $u_1^{(j+1)}\geq u_1^{(j)}$ 
and $u_2^{(j+1)}\geq u_2^{(j)}$ for some $j\in\N$, then 
\[\bald
(\gl+M) u_1^{(j+2)}+A_1(u_1^{(j+2)},u_2^{(j+2)})&\,=g_1+Mu_1^{(j+1)}\geq g_1+Mu_1^{(j)}, 
\\ (\gl+M)
 u_2^{(J+2)}+A_2(u_1^{(j+2)},u_2^{(j+2)})&\,=g_2+Mu_2^{(j+1)}\geq g_2+Mu_2^{(j)},
\eald
\]
and, by the comparison argument as above, we deduce that $u_1^{(j+2)}\geq u_1^{(j+1)}$ and 
$u_2^{(j+2)}\geq u_2^{(j+1)}$. By induction, we conclude that $u_1^{j+1}\geq u_1^{(j)}$ 
and $u_2^{(j+1)}\geq u_2^{(j)}$ for all $j\in\N$. 

We now know that the sequences $(u_1^{(j)})_{j\in\N}$ and $(u_2^{(j)})_{j\in\N}$ 
are convergent. Let $u_1$ and $u_2$ denote the respective limits, and we note that 
\[
(\gl+M) u_1+A_1(u_1,u_2)=Mu_1 \ \ \AND \ \ 
(\gl+M)u_2+A_2(u_1,u_2)=Mu_2,
\]
to conclude that $(u_1,u_2)$ is a solution of \erf{eq:3.2}. 
\eproof

We note that if $(X_\gl,Y_\gl)\in\R^2$ is the unique solution of \erf{eq:3} 
and if $H_1(0)=H_2(0)=0$, the pair of 
constant functions $u_1(x)=X_\gl$ and $u_2(x)=Y_\gl$ is a solution of \erf{eq:1}. 
As is well-known (see, for instance, \cite{IJ} and the references therein),
\erf{eq:1} has a unique (viscosity) solution, and hence, the pair $(X_\gl,Y_\gl)$ is the unique solution of \erf{eq:1}. 

The main result of this paper is the following two theorems. 

\begin{theorem} \label{thm1} For any $\gl>0$, let $(X_\gl,Y_\gl)\in\R^2$ be the solution of \erf{eq:3}. Then (i) the set of points $(X_\gl,Y_\gl)$, with $\gl>0$, 
is bounded in $\R^2$.  (ii) We have 
\[
\liminf_{\gl \to 0+}X_\gl\leq \fr{d}{2}<\limsup_{\gl \to 0+}X_\gl,
\ \ \AND \ \ 
\liminf_{\gl \to 0+}Y_\gl\leq -\fr{d}{2}< 
\limsup_{\gl \to 0+}Y_\gl.
\]
In particular, the family of the pairs $(X_\gl,Y_\gl)$ does not converge as 
$\gl \to 0+$.
\end{theorem}

As noted before the theorem, the following is an immediate consequence of 
Theorem \ref{thm1}. 

\begin{theorem} \label{cor1}Assume that $H_1(0)=H_2(0)=0$. 
%Let $B_i : \R^2\to\R$, with $i\in\I$, be the functions given by \erf{eq:3.1}
For any $\gl>0$, let $(u_{\gl,1},u_{\gl,2})$ be the (viscosity) solution of \erf{eq:1}. Then, 
the functions $u_{\gl,i}$ are constants, 
the family of the points $(u_{\gl,1},u_{\gl,2})$ 
in $\R^2$ is bounded, and it does not converge as $\gl\to 0+$.   
\end{theorem}

Notice that the convexity of $H_i$ in the above theorem is irrelevant, 
and, for example, one may take $H_i(p)=|p|^2$ for $i\in\I$, which are convex functions.  

We remark that a claim similar to Theorem \ref{cor1} is valid when one 
replaces $H_i(p)$ by degenerate elliptic operators $F_i(p,M)$ (see \cite{CIL} 
for an overview on the viscosity solution approach to 
fully nonlinear degenerate elliptic equations), where 
$M$ is the variable corresponding to the Hessian matrices of unknown functions.

In the next and final section, we give the proof of Theorem \ref{thm1}.

\bigskip

\section{Proof of Theorem \ref{thm1}}

This section is entirely devoted to the proof of Theorem \ref{thm1}.

For any  $\gl>0$, let $(X_\gl,Y_\gl)$ denote the solution of \erf{eq:3}.

\bproof[Proof of Theorem \ref{thm1}, (i)] We show first that $(d/2, -d/2)$ is a 
solution of 
\[
B_1(u_1,u_2)=0 \ \ \AND \ \ B_2(u_1,u_2)=0. 
\]
To see this, we observe that for any $\ga\in\cA$,
\[\bald
b_1(\ga,0,d/2,-d/2)&\,=(d+d^2)(\ga-\gamma),
\\ b_1(\ga,1,d/2,-d/2)&\,=-(1+d)(\ga-\gamma),
\\b_2(\ga,0,d/2,-d/2)&\,=-(d+d^2)(\ga-\gamma),
\\b_2(\ga,1,d/2,-d/2)&\,=(1+d)(\ga-\gamma),
\eald
\]
and hence, we get 
\[
B_1(d/2,-d/2)=0 \ \ \AND \ \ B_2(d/2,-d/2)=0.
\]
By \erf{eq:2}, we see that $B_i(u_1+r,u_2+r)=B_i(u_1,u_2)$ for all $i\in\I,\,(u_1,u_2)\in
\R^2$ and $r\in\R$. Hence, we have 
\[
B_i(d,0)=B_i(0,-d)=0 \ \ \FOR i\in\I,
\]
which shows that for any $\gl>0$, if $(v_1,v_2)=(d,0)$ (resp., $(v_1,v_2)=(0,-d)$) satisfies 
\[
\gl v_i+B_i(v_1,v_2)\geq 0 \ \ (\,\text{resp., }  \ \ \gl v_i+B_i(v_1,.v_2)\leq 0\, ) \ \ \FORALL i\in\I.
\] 
By the comparison assertion of Proposition \ref{prop1}, we deduce that
\[
0\leq X_\gl\leq d \ \ \AND \ \ -d\leq Y_\gl\leq 0 \ \ \FORALL \gl>0,
\]
which proves that the set of $(X_\gl,Y_\gl)$, with $\gl>0$, is bounded in $\R^2$. 
\eproof

\begin{lemma} \label{difference} We have 
\[
\lim_{\gl\to 0+}(X_\gl-Y_\gl)=d.
\]
\end{lemma} 

\bproof Set $Z_\gl=X_\gl-Y_\gl$ and 
\[
B(u)=\max_{\ga\in\cA}\min_{\gb\in\cB}
(c_{11}(\ga,\gb)u-L_1(\ga,\gb)) -\max_{\ga\in\cA}\min_{\gb\in\cB}(-c_{22}(\ga,\gb)u-L_2(\ga,\gb)).
\]
Noting that for any $(u_1,u_2)\in\R^2$,
\[\bald
B_1(u_1,u_2)&\,=\max_{\ga\in\cA}\min_{\gb\in\cB}\left(c_{11}(\ga,\gb)(u_1-u_2)-L_1(\ga,\gb)\right),
\\ B_2(u_1,u_2)&\,=\max_{\ga\in\cA}\min_{\gb\in\cB}\left(-c_{22}(\ga,\gb)(u_1-u_2)-L_1(\ga,\gb)\right),
\eald\]
we find that  $B_1(u_1,u_2)-B_2(u_1,u_2)=B(u_1-u_2)$ for all $(u_1,u_2)\in\R^2$, and 
\beq\label{diff.1}
\gl Z_\gl+B(Z_\gl)=0.
\eeq
Since $c_{11}(\ga,\gb)=c_{22}(\ga,\gb)=\ga+\gb-2\ga\gb
>(\ga-\gb)^2>0$ for all $(\ga,\gb)\in\cA\tim\cB$, 
the function $B$ is (strictly) increasing on $\R$. 

Thanks to the claim (i) of Theorem \ref{thm1}, the family $(Z_\gl)_{\gl>0}$ 
has a limit point $Z_0\in\R$ as $\gl\to 0+$. It follows from \erf{diff.1} 
and the continuity of $B$ that
\[
B(Z_0)=0. 
\]
Since $B$ is increasing, $Z_0$ is a unique zero of the function $B$, which implies
\[
\lim_{\gl\to 0+}Z_\gl=Z_0.
\]
On the other hand,  
the proof of (i) of Theorem \ref{thm1} shows that  $B_1(d/2,-d/2)=B_2(d/2,-d/2)=0$,
which implies that $B(d)=0$ and $Z_0=d$. Thus, we conclude that 
$\lim_{\gl\to 0+}Z_\gl=d$. 
\eproof

\bproof[Proof of Theorem \ref{thm1}, (ii)] By Lemma \ref{difference}, we have
\[
\liminf_{\gl\to 0+}Y_\gl=\liminf_{\gl\to 0+}X_\gl+\lim_{\gl \to 0+}(Y_\gl-X_\gl)
=\liminf_{\gl\to 0+}X_\gl -d. 
\] 
Similarly, we have
\[
\limsup_{\gl\to 0+}Y_\gl=\limsup_{\gl\to 0+}X_\gl -d.
\]
Hence, we only need to prove that 
\beq\label{eq:2.1}
\liminf_{\gl\to 0+}X_\gl \leq \fr d 2<\limsup_{\gl\to 0+}X_\gl.
\eeq

We fix $\ga_1,\ga_2\in(0,\,1)$ and $\gb_1,\gb_2\in\cB$, and consider the linear 
problem 
\beq\label{eq:2.2}\bcases
\gl X+c_{11}(\ga_1,\gb_1)X+c_{12}(\ga_1,\gb_1)Y-L_1(\ga_1,\gb_1)=0, &\\ 
\gl Y+c_{21}(\ga_2,\gb_2)X+c_{22}(\ga_2,\gb_2)Y-L_2(\ga_2,\gb_2)=0. &
\ecases
\eeq

When $(\gb_1,\gb_2)=(0,0)$, \erf{eq:2.2} reads
\[
\bcases
\gl X=-\ga_1 X+\ga_1 Y+d^2(1-\ga_1),& \\
\gl Y=\ga_2 X-\ga_2 Y-d^2(1-\ga_2).&
\ecases
\]
This yields
\[
\bald
X&\,=-\,\fr{d^2((\ga_1-1)\gl+\ga_1-\ga_2)}{\gl(\ga_1+\ga_2+\gl)},
\\ Y&\,=\fr{d^2((\ga_2-1)\gl-\ga_1 +\ga_2)}{\gl(\ga_1+\ga_2+\gl)}.
\eald
\]
When $(\gb_1,\gb_2)=(0,1)$, we have 
\[
\bcases
\gl X= -\ga_1 X+\ga_1 Y+d^2(1-\ga_1),& \\
\gl Y=(1-\ga_2) X-(1-\ga_2)Y-\ga_2, &
\ecases
\]
and 
\[
\bald
X&\,=\fr{-\ga_1\ga_2+d^2((-\ga_1+1)\gl +\ga_1\ga_2-\ga_1-\ga_2+1)}{\gl(\ga_1-\ga_2+\gl+1)},
\\ Y&\,=\fr{-\ga_2\gl -\ga_1\ga_2 +d^2(\ga_1\ga_2 -\ga_1-\ga_2+1)}{\gl(\ga_1-\ga_2+\gl+1)}.
\eald
\]
When $(\gb_1,\gb_2)=(1,0)$, 
\[
\bcases
\gl X=-(1-\ga_1)X+(1-\ga_1) Y+\ga_1, &\\
\gl)Y=\ga_2 X-\ga_2Y-d^2(1-\ga_2), &
\ecases
\]
and 
\begin{align*}
X&\,=\fr{\ga_1\gl +\ga_1\ga_2 +d^2(-\ga_1\ga_2+\ga_1+\ga_2-1)}{\gl(-\ga_1+\ga_2+\gl+1)},
\\ Y&\,=\fr{-\ga_1\ga_2 +d^2(( -\ga_2+1)\gl +\ga_1\ga_2 -\ga_1 -\ga_2+1)}{\gl(\ga_1-\ga_2-\gl-1)}.
\end{align*}
When $(\gb_1,\gb_2)=(1,1)$, we have
\[\bcases
\gl X=-(1-\ga_1)X+(1-\ga_1) Y+\ga_1, &\\
\gl Y=(1-\ga_2) X-(1-\ga_2)Y-\ga_2, &
\ecases
\]
and 
\[
\bald
X&\,=\fr{\ga_1\gl +\ga_1-\ga_2}{\gl(-\ga_1-\ga_2+\gl+2)},
\\ Y&\,=\fr{-\ga_1 +\ga_2\gl +\ga_2}{\gl(\ga_1+\ga_2-\gl-2)}.
\eald
\]
\bigskip

In what follows, for the solution $(X,Y)$ of \erf{eq:2.2}, we write 
\[
X=X(\gl,\ga_1,\ga_2,\gb_1,\gb_2),\quad Y=Y(\gl,\ga_1,\ga_2,\gb_1,\gb_2). 
\]

We set 
\[ 
\gth=\fr{4(d+1)}{d-1}, 
\]
and, for $n\in\N$, 
\[
\rho_n=4^{-n-K},\qquad
\gl_n=\gth 4^{-n-K}=\gth \rho_n. 
\]
We write  
\[
p(\gl)=\gamma+\gth^{-1}\gl \ \ \FOR \gl>0,
\]
and note that $p(\gl_n)=\gamma+\rho_n\in \cA$ and, since \
$d^2(\gamma^2-2\gamma+1)-\gamma^2=0$, 
\[\bald
d^2(p(\gl_n)^2-2p(\gl_n)+1)-p(\gl_n)^2
&\,=d^2(2(\gamma -1)\rho_n +\rho_n^2)-2\gamma\rho_n -\rho_n^2
\\&\,=-2d\rho_n +\rho_n^2(d^2-1).
\eald\]

We compute that
\begin{align*}
\lim_{n\to\infty}X(\gl_n,p(\gl_n),p(\gl_n),0,0)
&\,=-\lim_{n\to\infty}\fr{d^2\gl_n(p(\gl_n)-1)}{\gl_n(2p(\gl_n)+\gl_n)}
=-\fr{d^2(\gamma-1)}{2\gamma}=\fr{d}{2}, 
\\ 
\lim_{n\to\infty}
X(\gl_n,p(\gl_n),p(\gl_n),0,1) &\,=\lim_{n\to\infty}\fr{-p(\gl_n)^2+d^2((1-p(\gl_n))\gl_n+p(\gl_n)^2-2p(\gl_n)+1)}{\gl_n}
\\&\,=-d^2(\gamma-1)+\lim_{n\to\infty}\fr{d^2(p(\gl_n)^2-2p(\gl_n)+1)-p(\gl_n)^2}{\gth \rho_n}
\\&\, =\fr{d^2}{1+d}-2\gth^{-1}d= \fr{d}2, 
\\ \lim_{n\to \infty}X(\gl_n,p(\gl_n),p(\gl_n),1,0) 
&\, =\lim_{n\to\infty}\fr{p(\gl_n)\gl_n +p(\gl_n)^2 -d^2(p(\gl_n)^2-2p(\gl_n)+1)}{\gl_n}
\\&\,=\gamma -\lim_{n\to \infty}\fr{d^2(p(\gl_n)^2-2p(\gl_n)+1)-p(\gl_n)^2}{\gth\rho_n}
\\&\,= \gamma+2\gth^{-1}d=\fr d 2, 
\\\lim_{n\to\infty} X(\gl_n,p(\gl_n),p(\gl_n),1,1) 
&\,=\lim_{n\to\infty}\fr{p(\gl_n)\gl_n}{\gl_n(-2p(\gl_n)+\gl_n+2)}
=\fr{\gamma}{2(1-\gamma)}
=\fr{d}{2}.
\end{align*}

We set $(X_n,Y_n):=(X_{\gl_n},Y_{\gl_n})$, so that 
\[
\gl_n X_n+B_1(X_n,Y_n)=0 \ \ \AND \ \ \gl_n Y_n+B_2(X_n,Y_n)=0.
\]
We have for any $\ga\in\cA$, 
\[
%\bald&
\gl_n X_n+\min_{\gb\in\cB}b_1(\ga,\gb,X_n,Y_n)\leq 0,
%\\&
\ \ \AND \ \  \gl_n Y_n+\min_{\gb\in\cB}b_2(\ga,\gb,X_n,Y_n)\leq 0,
%\eald
\]
and, in particular,
\[
\gl_n X_n+\min_{\gb\in\cB}b_1(p(\gl_n),\gb,X_n,Y_n)\leq 0 
 \ \ \AND \ \ \gl_n Y_n+\min_{\gb\in\cB}b_2(p(\gl_n),\gb,X_n,Y_n)\leq 0.
\]
We select $\gb_{n,1},\gb_{n,2}\in\cB$ so that 
\[
\gl_n X_n+b_1(p(\gl_n),\gb_{n,1},X_n,Y_n)\leq 0
\ \ \AND \ \ \gl_n Y_n+b_2(p(\gl_n),\gb_{n,2},X_n,Y_n)\leq 0.
\] 
The comparison assertion of Proposition \ref{prop1}, with 
$\cA_1=\cA_2=\{p(\gl_n)\}$ and $\cB_i=\{\gb_{n,i}\}$, yields
\[
X_n\leq X(\gl_n,p(\gl_n),p(\gl_n),\gb_{n,1},\gb_{n,2})\ \ \AND \ \ 
Y_n\leq Y(\gl_n,p(\gl_n),p(\gl_n),\gb_{n,1},\gb_{n,2}).
\]
Combining all together, we find that 
\[
\limsup_{n\to \infty}X_n\leq \fr{d}{2}.
\]
In particular, we have 
\beq\label{eq:2.3}
\liminf_{\gl\to 0+}X_\gl\leq \fr{d}{2}.
\eeq

To proceed the proof, we check the monotonicity of $X(\gl,\ga_1,\ga_2,\gb_1,\gb_2)$ as a function of $\ga_1,\,\ga_2$. 
We use the notation: for $Z=X$ or $Z=Y$, 
\begin{align*}
&\pl_1 Z(\gl,\ga_1,\ga_2,\gb_1,\gb_2)=\pl_t Z(\gl,t,\ga_2,\gb_1,\gb_2)\Big|_{t=\ga_1}, 
\\&\pl_2 Z(\gl,\ga_1,\ga_2,\gb_1,\gb_2)=\pl_t Z(\gl,\ga_1,t,\gb_1,\gb_2)\Big|_{t=\ga_2}.
\end{align*}
By simple computation, we obtain 
\beq \label{eq:2.4}\bald
\pl_{1} X(\gl,\ga_1,\ga_2,0,0)&\,=
-\fr{d^2(\gl+2)(\gl+\ga_2)}{\gl(\ga_1+\ga_2+\gl)^2}<0,
\\ \pl_{2} X(\gl,\ga_1,\ga_2,0,0)&\,=\fr{d^2\ga_1(\gl+2)}{\gl(\ga_1+\ga_2+\gl)^2}>0, 
\\ \pl_{1} X(\gl,\ga_1,\ga_2,0,1)&\,=
-\fr{(\gl+1-\ga_2)(d^2(\gl+2-\ga_2)+\ga_2)}{\gl(\ga_1-\ga_2+\gl+1)^2}<0,
\\ \pl_{2} X(\gl,\ga_1,\ga_2,0,1)&\,=
-\fr{\ga_1(\gl+d^2(1-\ga_1)+\ga_1+1)}{\gl(\ga_1-\ga_2+\gl+1)^2}<0.
\\ \pl_{1} X(\gl,\ga_1,\ga_2,1,0)&\,=
\fr{(\gl+\ga_2)(\gl+d^2(1-\ga_2)+\ga_2+d^2+1)}{\gl(\ga_1-\ga_2+\gl+1)^2}>0,
\\ \pl_{2} X(\gl,\ga_1,\ga_2,1,0)&\,
=\fr{(1-\ga_1)(d^2(\gl-\ga_1+2)+\ga_1)}{\gl(\ga_1-\ga_2+\gl+1)^2}>0.
\\ \pl_{1} X(\gl,\ga_1,\ga_2,1,1)&\,=
\fr{(\gl+2)(\gl-\ga_2+1)}{\gl(-\ga_1-\ga_2+\gl+2)^2}>0,
\\ \pl_{2} X(\gl,\ga_1,\ga_2,1,1)&\,=
-\fr{(1-\ga_1)(\gl+2)}{\gl(-\ga_1-\ga_2+\gl+2)^2}<0.
\eald
\eeq

\allowdisplaybreaks[4]

Now, set
\[
\tau:=\fr{10(d+1)}{d-1},
\qquad\rho_n=4^{-n-K},
\qquad \mu_n=2\tau \rho_n,
\]
and write $\,q(\gl)=\gamma+\tau^{-1}\gl$ for $\gl>0$.
Note that
\begin{gather*}
q(\mu_n/2)=\gamma+\rho_n,
\qquad q(2\mu_n)=\gamma+4\rho_n, \qquad\qquad 
\\ d^2(q(\mu_n/2)^2-2q(\mu_n/2)+1)-q(\mu_n/2)^2=
-2d\rho_n+(d^2-1)\rho_n^2, 
\\ d^2(q(2\mu_n)^2-2q(2\mu_n)^2+1)-q(2\mu_n)^2=
-8d\rho_n+16(d^2-1)\rho_n^2, 
\\ d^2(q(\mu_n/2)q(2\mu_n)-q(\mu_n/2)-q(2\mu_n)+1)-q(\mu_n/2)q(2\mu_n)=
-5d\rho_n+4(d^2-1)\rho_n^2.
\end{gather*}

We compute by using the above equalities that 
\begin{align*}
\lim_{n\to\infty}&X(\mu_n,q(\mu_n/2),q(2\mu_n),0,0)=
\\&= -\lim_{n\to\infty}\fr{d^2(q(\mu_n/2)-1)}{q(\mu_n/2)+q(2\mu_n)} 
-\lim_{n\to\infty}\fr{d^2(q(\mu_n/2)-q(2\mu_n))}{\mu_n(q(\mu_n/2)+q(2\mu_n))}
\\&=\fr{d^2(1-\gamma)}{2\gamma}+\fr{3d^2}{2\tau \gamma}
=\fr d{2}+\fr{3d(d-1)}{20}>\fr d{2}, 
%%%%%
\\ \lim_{n\to\infty}&X(\mu_n,q(\mu_n/2),q(\mu_n/2),0,1)
\\& =\lim_{n\to\infty}d^2(-q(\mu_n/2)+1)
\\&\quad +\lim_{n\to\infty}\fr{-q(\mu_n/2)^2+d^2(q(\mu_n/2)^2-2q(\mu_n/2)+1)}{\mu_n}
\\&=d^2(1-\gamma) 
-\fr{2d}{\tau}
=\fr{d}{2}+\fr{3d(d-1)}{10(d+1)}>\fr{d}2,
%%%%%%%%%%%
\\ \lim_{n\to\infty}&X(\mu_n,q(2\mu_n),q(2\mu_n),1,0)
\\&= \lim_{n\to\infty}q(2\mu_n)
-\lim_{n\to\infty}\fr{-q(2\mu_n)^2+d^2(q(2\mu_n)^2-2q(2\mu_n)+1)}{\mu_n}
\\&=\gamma +\fr{8d}{\tau}=\fr d2 +\fr{3d(d-1)}{10(d+1)}>\fr{d}{2}, 
%%%%%%%%%%%%%%%
\\ \lim_{n\to\infty}&X(\mu_n,q(2\mu_n),q(\mu_n/2),1,1)
\\&= \lim_{n\to\infty}\fr{q(2\mu_n)}
{-q(2\mu_n)-q(\mu_n/2)+2}+
\lim_{n\to\infty}\fr{q(2\mu_n)-q(\mu_n/2)}
{\mu_n(-q(2\mu_n)-q(\mu_n/2)+2)}
\\&=\fr{\gamma}{2(1-\gamma)}+\fr {3}{2\tau(1-\gamma)}
=\fr d2 +\fr{3(d-1)}{20}>\fr d{2}.
\end{align*}

We set 
\[\bald
x_n&\,=\min\{X(\mu_n,q(\mu_n/2),q(2\mu_n),0,0),X(\mu_n,q(\mu_n/2),q(\mu_n/2),0,1),
\\&\,\qquad \qquad X(\mu_n,q(2\mu_n),q(2\mu_n),1,0),X(\mu_n,q(2\mu_n),q(\mu_n/2),1,1)\}.
\eald\]

\def\hX{\widehat X} \def\hY{\widehat Y} 
\def\hga{\hat\ga}\def\hgb{\hat\gb}
For any $n\in\N$, set $(\hX_n,\hY_n):=(X_{\mu_n},Y_{\mu_n})$ and note 
\[
\mu_n \hX_n+B_1(\hX_n,\hY_n)=0 \ \ \AND \ \ 
\mu_n \hY_n+B_2(\hX_n,\hY_n)=0.
\]
We select $\hga_{n,1},\hga_{n,2}\in\cA$ and $\hgb_{n,1},\hgb_{n,2}\in\cB$ so that 
\[
\bcases
\mu_n\hX_n+b_1(\hga_{n,1},\hgb_{n,1},\hX_n,\hY_n)=0,&\\
\mu_n\hY_n+b_2(\hga_{n,2},\hgb_{n,2},\hX_n,\hY_n)=0.
\ecases
\]
Accordingly, we have 
\[
\hX_n=X(\mu_n,\hga_{n,1},\hga_{n,2},\hgb_{n,1},\hgb_{n,2}) \ \ \AND \ \ 
\hY_n=X(\mu_n,\hga_{n,1},\hga_{n,2},\hgb_{n,1},\hgb_{n,2}). 
\]

Note that $q(\mu_n/2),q(2\mu_n)\in\cA$ and that 
if $t\in (q(\mu_n/2),\,q(2\mu_n)$, then $t\not\in\cA$.
The monotonicity of $X$ (see \erf{eq:2.4}) shows that for any $\ga_1,\ga_2\in\cA$, 
if $\ga_1\geq q(\mu_n)$ and $\ga_2\geq q(\mu_n)$, then  
\[
X(\mu_n,\ga_1,\ga_2,1,0)\geq X(\mu_n,q(2\mu_n),q(2\mu_n),1,0)\geq x_n,
\]
if $\ga_1> q(\mu_n)$ and $\ga_2\leq q(\mu_n)$, then 
\[
X(\mu_n,\ga_1,\ga_2,1,1)\geq X(\mu_n,q(2\mu_n),q(\mu_n/2),1,1)\geq x_n,
\]
if $\ga_1\leq q(\mu_n)$ and $\ga_2\geq q(\mu_n)$, then 
\[
X(\mu_n,\ga_1,\ga_2,0,0)\geq X(\mu_n,q(\mu_n/2),q(2\mu_n),0,0)\geq x_n,
\]
and, if $\ga_1\leq q(\mu_n/2)$ and $\ga_2\leq q(\mu_n/2)$, then 
\[
X(\mu_n,\ga_1,\ga_2,0,1)\geq X(\gl,q(\mu_n/2),q(\mu_n/2),0,1)\geq x_n.
\]
Hence, we deduce that 
\[
\hX_n=X(\mu_n,\hga_{n,1},\hga_{n,2},\hgb_{n,1},\hgb_{n,2})\geq x_n. 
\]
From these, we conclude that 
\[
\liminf_{n\to\infty} \hX_n\geq \liminf_{n\to\infty}x_n>\fr d{2},
\]
which ensures that
\[
\limsup_{\gl\to 0+}X_\gl>\fr d{2}. 
\]
The proof is now complete.
\eproof 

\begin{proposition} \label{prop.=d/2} We have 
\[
\liminf_{\gl\to 0+}X_\gl=\fr d2 \ \ \AND \ \ \liminf_{\gl\to 0+}Y_\gl=-\fr d2. 
\]
\end{proposition}

\bproof  As in the proof of \erf{eq:2.3}, we set 
\[ 
\gth=\fr{4(d+1)}{d-1},  \ \ \AND \ \  
p(\rho)=\gamma+\rho \ \ \FOR \rho>0. 
\]

Let $X(\ga_1,\ga_2,\gb_1,\gb_2)$ be the function 
as in the proof of Theorem \ref{thm1}, (ii) for $\ga_1,\ga_2\in (0,1)$
and $\gb_1,\gb_2\in\cB$.
Noting that $p(\rho)=\gamma+\rho\in (0,1)$ and that
\[\bald
d^2(p(\rho)^2-2p(\rho)+1)-p(\rho)^2
&\,=d^2(2(\gamma -1)\rho +\rho^2)-2\gamma\rho -\rho^2
\\&\,=-2d\rho +\rho^2(d^2-1),
\eald\]
we compute that 
\begin{align*}
\lim_{\rho \to 0+}X(\gth\rho,p(\rho),p(\rho),0,0)
&\,=-\lim_{\rho\to0+}\fr{d^2 \gth\rho(p(\rho)-1)}{\gth\rho (2p(\rho)+\gth\rho)}
=-\fr{d^2(\gamma-1)}{2\gamma}=\fr{d}{2}, 
\\ 
\lim_{\rho\to 0+}
X(\gth\rho,p(\rho),p(\rho),0,1) &\,=\lim_{\rho \to 0+}\fr{-p(\rho)^2+d^2(
(1-p(\rho))\gth\rho+p(\rho)^2-2p(\rho)+1)}{\gth\rho}
\\&\,=d^2(1-\gamma)-\fr{2d}{\gth}= \fr{d}2, 
\\ \lim_{\rho\to 0+}X(\gth\rho,p(\rho),p(\rho),1,0) 
&\, =\lim_{\rho \to 0+}\fr{p(\rho)\gth
\rho +p(\rho)^2 -d^2(p(\rho)^2-2p(\rho)+1)}{\gth\rho}
\\&\,=\gamma -\lim_{\rho\to 0+}\fr{d^2(p(\rho)^2-2p(\rho)+1)-p(\rho)^2}{\gth\rho}
\\&\,= \gamma+2\gth^{-1}d=\fr d 2, 
\\\lim_{\rho\to 0+} X(\gth\rho,p(\rho),p(\rho),1,1) 
&\,=\lim_{\rho \to 0+}\fr{p(\rho)\gth\rho}{\gth\rho(-2p(\rho)+\gth\rho+2)}
=\fr{\gamma}{2(1-\gamma)}
=\fr{d}{2}.
\end{align*}
Thus, setting 
\[\bald
x(\rho)=\min\{&X(\gth\rho,p(\rho),p(\rho),1,0),
X(\gth\rho,p(\gl),p(\gl),1,1), \\
&X(\gth\rho,p(\rho),p(\rho),0,0),
X(\gth\rho,p(\rho),p(\rho),0,1)
\},
\eald \]
we have 
\[
\lim_{\rho\to 0+}x(\rho)=\fr d2. 
\]

For any $\gl>0$,
we select $\ga_{\gl,1},\ga_{\gl,2}\in\cA$ and $\gb_{\gl,1},\gb_{\gl,2}\in\cB$ so that 
\[
\bcases
\gl X_{\gl}+b_1(\ga_{\gl,1},\gb_{\gl,1},X_\gl,Y_\gl)=0,&\\
\gl Y_\gl+b_2(\ga_{\gl,2},\gb_{\gl,2},X_\gl,Y_\gl)=0.
\ecases
\]
Accordingly, we have 
\[
X_\gl=X(\gl,\ga_{\gl,1},\ga_{\gl,2},\gb_{\gl,1},\gb_{\gl,2}) \ \ \AND \ \ 
Y_\gl=X(\gl,\ga_{\gl,1},\ga_{\gl,2},\gb_{\gl,1},\gb_{\gl,2}). 
\]

The monotonicity of $X$ (see \erf{eq:2.4}) shows that for any $\ga_1,\ga_2\in\cA$, 
if $\ga_1\geq p(\rho)$ and $\ga_2\geq p(\rho)$, then  
\[
X(\gth\rho,\ga_1,\ga_2,1,0)\geq X(\gth\rho,p(\rho),p(\rho),1,0)\geq x(\rho),
\]
if $\ga_1> p(\rho)$ and $\ga_2\leq p(\rho)$, then 
\[
X(\gth\rho,\ga_1,\ga_2,1,1)\geq X(\gth\rho,p(\gl),p(\gl),1,1)\geq x(\rho),
\]
if $\ga_1\leq p(\rho)$ and $\ga_2\geq p(\rho)$, then 
\[
X(\gth\rho,\ga_1,\ga_2,0,0)\geq X(\gth\rho,p(\rho),p(\rho),0,0)\geq x(\rho),
\]
and, if $\ga_1\leq p(\rho)$ and $\ga_2\leq p(\rho)$, then 
\[
X(\gth\rho,\ga_1,\ga_2,0,1)\geq X(\gth\rho,p(\rho),p(\rho),0,1)\geq x(\rho).
\]
Hence, 
we deduce that 
\[
X_{\gth\rho}=X(\gth\rho,\ga_{\rho,1},\ga_{\rho,2},\gb_{\rho,1},\gb_{\rho,2})\geq x(\rho), 
\]
and conclude that 
\[
\liminf_{\gl\to 0+}X_\gl \geq \lim_{\rho\to 0+}x(\rho)=\fr d{2},
\]
which finishes the proof.
\eproof

\section*{Acknowledgments}
%The authors would like to thank the anonymous referees for their careful reading 
%of and critical and useful comments on the original version of this paper, which 
%have helped significantly to improve the presentation.  
The author was supported in part by the JSPS Grants KAKENHI  No. 16H03948, 
No. 20K03688, No. 20H01817 and No. 21H00717. 
He thanks Wolfram Alpha for helping him checking  his computation.

\bye